\documentclass{amsart}
\usepackage{graphicx}
\theoremstyle{definition}

\theoremstyle{remark}

\numberwithin{equation}{section}

\begin{document}

\title{A note on K$\ddot{A}$hler manifolds with almost nonnegative bisectional curvature}%
\author{Hong Huang}%
\address{School of Mathematical Sciences,Key Laboratory of Mathematics and Complex Systems,
Beijing Normal University,Beijing 100875, P. R. China}%
\email{hhuang@bnu.edu.cn}%

\thanks{Partially supported by NSFC no.10671018.}%
\subjclass{53C44}%

\keywords{K$\ddot{a}$hler manifolds, almost nonnegative bisectional curvature, Ricci flow}%

\begin{abstract}
In this note we prove the following result: There is a positive
constant $\varepsilon(n,\Lambda)$ such that if $M^n$ is a simply
connected compact K$\ddot{a}$hler manifold with sectional curvature
bounded from above by $\Lambda$, diameter bounded from above by 1,
and with holomorphic bisectional curvature $H \geq
-\varepsilon(n,\Lambda)$, then $M^n$ is diffeomorphic to the product
$M_1\times ... \times M_k$, where each $M_i$ is either a complex
projective space or an irreducible K$\ddot{a}$hler-Hermitian
symmetric space of rank $\geq 2$. This resolves  a conjecture of F.
Fang under the additional upper bound restrictions on sectional
curvature and diameter.
\end{abstract} \maketitle

\section {Introduction}

In a recent paper [PT] Petersen and Tao give a classification of
simply connected compact manifold with almost quarter-pinched
curvature using the Ricci flow method. Inspired by their work, in
this short note we will try to classify simply connected compact
K$\ddot{a}$hler manifolds with almost nonnegative bisectional
curvature.

First recall (cf. Fang [F]) that a compact complex manifold $M$ has
almost nonnegative bisectional curvature if given any $\varepsilon
>0$, there is a Hermitian metric $g$ on $M$ with bisectional
curvature $H$ satisfying $H \cdot diam(M_g)^2 \geq -\varepsilon$.

 Our main result is the following

\hspace *{0.4cm}

  {\bf Theorem } There is a positive constant
$\varepsilon(n,\Lambda)$ such that if $M^n$ is a simply connected
compact K$\ddot{a}$hler manifold with sectional curvature bounded
from above by $\Lambda$, diameter bounded from above by 1, and with
holomorphic bisectional curvature $H \geq -\varepsilon(n,\Lambda)$,
then $M^n$ is diffeomorphic to the product $M_1\times ... \times
M_k$, where each $M_i$ is either a complex projective space or an
irreducible K$\ddot{a}$hler-Hermitian symmetric space of rank $\geq
2$.

\hspace *{0.4cm}

 This improves earlier work of Fang [F].

 When we impose an additional condition  that the second Betti number
 $b_2(M^n)=1$,  we have
 the following

\hspace *{0.4cm}

 {\bf Corollary} There is a positive
constant $\varepsilon(n,\Lambda)$ such that if $M^n$ is a simply
connected compact K$\ddot{a}$hler manifold with sectional curvature
bounded from above by $\Lambda$, diameter bounded from above by 1,
holomorphic bisectional curvature $H \geq -\varepsilon(n,\Lambda)$,
and with the second Betti number $b_2(M^n)=1$, then $M^n$ is
diffeomorphic to either the complex projective space or an
irreducible K$\ddot{a}$hler-Hermitian symmetric space of rank $\geq
2$.

\hspace *{0.4cm}

 This resolves a conjecture of Fang [F] under the
additional upper bound restrictions on sectional curvature and
diameter.

 As in [PT] we also use the Ricci flow method in the proof of our theorem,
 which is given in the following section.

\section {Proof of Theorem}

 We prove by contradiction. Suppose the result is not true. Then
 we can find a sequence of simply
connected compact K$\ddot{a}$hler manifolds $(M_i,g_i)$ of complex
dimension $n$ satisfying the following condition:

i) the bisectional curvature $H_i$ of $M_i$ satisfies $H_i\geq
-\frac{1}{i}$, the sectional curvature $K_i$ of $M_i$ satisfies
$K_i\leq \Lambda$, and the diameter $diam(M_i)\leq 1$;

ii) none of  $M_i$ is  diffeomorphic to a product where each factor
is a complex projective space or an irreducible
K$\ddot{a}$hler-Hermitian symmetric space of rank $\geq 2$.

Note that the sectional curvature of $(M_i,g_i)$ is uniformly
bounded from below by a constant $\lambda$. Then by Hamilton's
work [H1] we can run the K$\ddot{a}$hler-Ricci flow $(M_i,g_i(t))$
for a fixed amount of time $[0,t_0](t_0>0)$ with initial data
$(M_i,g_i)$, with the sectional curvature of $(M_i,g_i(t))$
uniformly bounded on this time interval.

Also note that by a result of Fang and Rong [FR], the injectivity
radius of $(M_i,g_i)$ is uniformly bounded away from zero. Then by
 Hamilton's compactness theorem  for the Ricci flow [H3], $(M_i,g_i(t))$
  sub-converge to a K$\ddot{a}$hler Ricci flow $(M,g(t)),t\in [0,t_0]$.
The convergence is in the $C^{1,\alpha}$($\alpha <1$) sense at
$t=0$, and is in the $C^\infty$ sense for $0<t\leq t_0$.

Then by virtue of Hamilton's work (see the proof of [H2, Theorem
4.3]; compare with [PT]), the bisectional curvature of
$(M_i,g_i(t))$ is bounded from below by $-\frac{exp(Ct)}{i}$, where
$C$ is a constant depending only on the curvature bounds for
$g_i(0)$. (Note that here we have used implicitly  the work of Bando
[B] (in complex dimension 3) and Mok [M](in any dimension) on
preserving   nonnegative bisectional curvature condition  under the
K$\ddot{a}$hler-Ricci flow.)  So the limit $(M,g(t))$ has
nonnegative bisectional curvature for $t>0$. Then by Mok's theorem
[M] $M$ is diffeomorphic to a product, where each factor is a
complex projective space or an irreducible K$\ddot{a}$hler-Hermitian
symmetric space of rank $\geq 2$. But for $i$ sufficiently large,
$M_i$ is diffeomorphic to $M$, so $M_i$ is also diffeomorphic to
such a product, a contradiction, and we are done.

\bibliographystyle{amsplain}

\hspace *{0.4cm}

{\bf Reference}

\bibliography{1}[B] S. Bando, On the classification of three-dimensional compact Kaehler manifolds
of nonnegative bisectional curvature, J. Diff. Geom. 19 (1984),
no.2, 283-297.

\bibliography{2}[F] F. Fang, K$\ddot{a}$hler manifolds with almost
non-negative bisectional curvature, Asian J. Math. 6 (2002), no.3,
385-398.

\bibliography{3}[FR] F. Fang, X. Rong,  Fixed point free circle actions and finiteness
theorems, Comm. Contemp. Math. 2(2000), 75-96.

\bibliography{4}[H1] R. Hamilton, Three-manifolds with positive
Ricci curvature, J. Diff. Geom. 17(1982), 255-306.

\bibliography{5}[H2] R. Hamilton, Four-manifolds with positive
curvature operator, J. Diff. Geom. 24(1986), 153-179.

\bibliography{6}[H3] R. Hamilton, A compactness property for
solutions of the Ricci flow, Amer. J. Math. 117 (1995), no.3,
545-572.

\bibliography{7}[M] N. Mok, The uniformization theorem for compact
K$\ddot{a}$hler manifolds of  non-negative holomorphic bisectional
curvature, J. Diff. Geom.27 (1988), no.2,179-214.

\bibliography{8}[PT] P. Petersen, T. Tao, Classifiacation of almost quarter-pinched manifolds,
arXiv:0807.1900v1.

\end{document}